\documentclass[a4paper]{article}
\usepackage[T1]{fontenc}
\usepackage[utf8]{inputenc}
\usepackage{amsmath}
\usepackage{amsfonts}
\usepackage{amssymb}
\usepackage{graphicx}
\usepackage{tikz}
\usepackage{mathtools}
\usepackage[top=3cm, bottom=3cm, left=3cm, right=3cm]{geometry}
\usepackage{calrsfs}
\usepackage{bbm}
\usepackage{algorithm}
\usepackage[noend]{algpseudocode}
\usepackage{float}
\usepackage{enumerate}
\DeclareMathAlphabet{\pazocal}{OMS}{zplm}{m}{n}
\usepackage{lipsum}

\newcommand\blfootnote[1]{%
  \begingroup
  \renewcommand\thefootnote{}\footnote{#1}%
  \addtocounter{footnote}{-1}%
  \endgroup
}

\makeatletter
\newcommand*\bigcdot{\mathpalette\bigcdot@{.5}}
\newcommand*\bigcdot@[2]{\mathbin{\vcenter{\hbox{\scalebox{#2}{$\m@th#1\bullet$}}}}}
\makeatother

\newcommand\mydef{\stackrel{\mathclap{\normalfont\mbox{def}}}{=}}

\newtheorem{theorem}{\indent Theorem}[section]
\newtheorem{proposition}[theorem]{\indent Proposition}
\newtheorem{remark}[theorem]{\indent Remark}
\newtheorem{lemma}[theorem]{\indent Lemma}

\newtheorem{definition-theorem}[theorem]{\indent Definition-Theorem}
\newtheorem{corollary}[theorem]{\indent Corollary}

\newenvironment{proof}{\paragraph{Proof:}}{\hfill$\square$}

\def \N{\mathbb{N}}
\def \R{\mathbb{R}}
\def \Z{\mathbb{Z}}

\def \G{\mathcal{G}}
\def \F{\mathcal{F}}
\def \L{\mathcal{L}}
\def \I{\mathcal{I}}

\def \P{\mathbb{P}}
\newcommand{\E}{\mathbb{E}}

\title{\bf A Result of Metastability for an Infinite System of Spiking Neurons}

\author{Morgan André \\ \textit{Instituto de Matemática e Estatística,} \\ \textit{Universidade de São Paulo.}}

\begin{document}

\maketitle

\blfootnote{© 2020. This manuscript version is made available under the CC-BY-NC-ND 4.0 license}

\begin{abstract}
In 2018, Ferrari et al. wrote a paper called "Phase Transition for Infinite Systems of Spiking Neurons" in which they introduced a continuous time stochastic model of interacting neurons. This model consists in a countable number of neurons, each of them having an integer-valued membrane potential, which value determine the rate at which the neuron spikes. This model has also a parameter $\gamma$, corresponding to the rate of the leak times of the neurons, that is, the times at which the membrane potential of a given neuron is spontaneously reset to its resting value (which is $0$ by convention). As its title says, it was proven in this previous article that this model presents a phase transition phenomenon with respect to $\gamma$. Here we prove that this model also exhibits a metastable behavior. By this we mean that if $\gamma$ is small enough, then the re-normalized time of extinction of a finite version of this system converges toward an exponential random variable of mean $1$ as the number of neurons goes to infinity.
\end{abstract}

\vspace{0.4 cm}


\section{Introduction}

\vspace{0.4 cm}
In the present paper we consider an infinite system of spiking neurons introduced by Ferrari et al. in \cite{ferrari}. Informally this model can be described as follows: we have a countable set of neurons $I$, and to each neuron $i \in I$ is associated a set of \textit{presynaptic neurons} $\mathbb{V}_i$ and a process $(X_i(t))_{t \geq 0}$ which represents the \textit{membrane potential} of neuron $i$. Moreover, we associate to each neuron a Poisson process $(N^{\dagger}_i(t))_{t \geq 0}$ of some parameter $\gamma$, representing the \textit{leak times}. At any of these leak times the membrane potential of the neuron concerned is reset to $0$. Another point process $(N_i(t))_{t \geq 0}$ representing the \textit{spiking times} is also associated to each neuron, which rate at time $t$ is given by $\phi_i (X_i(t))$, where $\phi_i$ is some rate function that needs to be specified (typically a hard threshold, a linear function or a sigmoid function). Whenever a neuron spikes its membrane potential is also reset to $0$ and the membrane potential of all of its post-synaptic neurons is increased by one (i.e. the neurons of the set $\{ j : i \in  \mathbb{V}_j\}$).  We refer to section 2 of \cite{ferrari} for a more formal description of the model.

\vspace{0.4 cm}

This model can be seen as a variant in continuous time of the model introduced by A. Galves and E. Löcherbach in \cite{galves}, sometimes called the Galves-Löcherbach model, or simply GL model. Various other variants of this model have been discussed in the literature, and we refer to \cite{galves2016} for a review. Contrary to classical leaky integrate-and-fire models this model presents an inherent stochasticity. The leakage effect - that is to say the natural diffusion of ions which occurs through the membrane, resulting in a decrease of membrane potential - is often modeled as a continuous decrease via some negative term in a differential equation. In our model a different approach is considered as the leakage only occurs at discrete times, and it is abrupt as it immediately resets the membrane potential to $0$. This is mathematically convenient as it allows us to work only with point processes. The idea of using point processes to describe biological neural networks is not new and to the best of our knowledge it likely started with D. Brillinger in \cite{brillinger}, as well as with A.G. Hawkes in \cite{hawkes}. Nonetheless most of the studies related to Hawkes processes consider only finite systems. The system we consider here is infinite, as well as all the variations of the model introduced in \cite{galves} studied up to now. We believe that it is an interesting and reasonable approach as an actual biological neural system have a huge number of components (of the order of $10^{11}$).

\vspace{0.4 cm}

Here, as it was done in \cite{ferrari}, we consider a specific instantiation of this model where $I = \Z$ (where $\Z$ denotes the set of all integers), $\mathbb{V}_i = \{i-1,i+1\}$ and $\phi_i(x) = \mathbbm{1}_{x > 0}$ for all $i \in \Z$. In this paradigm, at any time $t\geq 0$, a neuron $i$ is said to be \textit{active} when its membrane potential $X_i(t)$ is strictly greater than $0$ and \textit{quiescent} when it's equal to $0$. The reason for this terminology is that for our choice of $\phi$ whenever the membrane potential of a given neuron is positive there is a positive probability that the neuron will spike in the the near future while this probability is null when the membrane potential is null. It was shown in \cite{ferrari} that this instantiation of the model presents a phase transition. More precisely the following theorem was proven.

\begin{theorem} \label{thm:phasetransition}
Suppose that for any $i \in \Z$ we have $X_i(0) \geq 1$. There exists a critical value $\gamma_c$ for the parameter $\gamma$, with $0 < \gamma_c < \infty$, such that for any $i \in \Z$

$$\P \Big( N_i([0,\infty[) \text{ } < \infty \Big) = 1 \text{ if } \gamma > \gamma_c$$

and

$$ \P \Big( N_i([0,\infty[) \text{ } = \infty \Big) > 0 \text{ if } \gamma < \gamma_c.$$

\end{theorem}

\vspace{0.4 cm}

In words the system continues spiking infinitely often with positive probability when $\gamma$ is small enough, and it stops spiking once for all after some time (at least locally) when $\gamma$ is big enough. The last sentence can be rephrased saying that for $\gamma > \gamma_c$ all neurons become quiescent as the time goes to infinity while for $\gamma < \gamma_c$ every neuron stays infinitely often active with positive probability.

\vspace{0.4 cm}

Our main result is that in the sub-critical regime (or at least in a sub-region of the sub-critical regime) the re-normalized time of extinction of the finite version of this system converges to an exponential random variable of unit mean when the number of neurons goes to infinity. More formally, we consider the process defined on a finite compact window $I_N \text{ } \mydef \text{ } \Z \cap [-N,N]$ for some $N \in \N$, instead of the whole lattice $\Z$. We define the extinction time of this process $$\tau_N \text{ } \mydef \text{ } \inf \Big\{t \geq 0 : X^N_i(t) = 0 \text { for all } i \in \Z \cap [-N,N]\Big\},$$ where $X^N_i(t)$ denotes the membrane potential of neuron $i$ at time $t$ in the finite system defined on $I_N$. We show that there exists $\gamma'_c$ such that if $0 < \gamma < \gamma'_c$, then we have $$\frac{\tau_N}{\E (\tau_N)} \overset{\mathcal{D}}{\underset{N \rightarrow \infty}{\longrightarrow}} \mathcal{E} (1),$$ where $\E$ denotes the mathematical expectation, the superscript $\mathcal{D}$ denotes a convergence in distribution and $\mathcal{E} (1)$ denotes an exponential random variable of mean $1$.

\vspace{0.4 cm}

In order to do this we consider an auxiliary process, namely the system of \textit{spiking rates} of the neurons. In our model each neuron in the one-dimensional lattice $\Z$ has only two possible spiking rates which are $1$ and $0$, depending on whether the neuron is active or quiescent. Any active neuron can be affected by two different effects at random exponential times: the occurrence of a \textit{spike}, and the \textit{leakage} effect. When a neuron is active the spikes occur as the atoms of a Poisson process of parameter 1, and when a spike occurs the neuron instantaneously becomes quiescent ($1 \rightarrow 0$) and his two post-synaptic neurons (which are his immediate neighbours on the right and on the left on the lattice) instantaneously become active if they weren't already ($0 \rightarrow 1$). Furthermore, an active neuron can becomes quiescent (without transmitting his activity to any neighbours) if it is affected by one of the leakage events, which, for a given neuron, occur as the atoms of a Poisson process of parameter $\gamma$.
 
\vspace{0.4 cm}

This spiking rates process can be seen as an interacting particle system with one single parameter $\gamma$. It has the important feature of being \textit{additive} (see \cite{harris78}) and for this reason it has a \textit{dual process} (see \cite{bertein}), which will be crucial to our purpose. The extinction time of the original model and the extinction time of the auxiliary process both correspond to the first time when all neurons are quiescent, thus these two are trivially equal. Moreover we notice that the auxiliary process is actually a continuous time Markov chain with a finite state space and with an absorbent state (the state where all neurons are quiescent), so that from elementary results on Markov chains it is clear that it will die out almost surely for any value of $\gamma$, which means that $\tau_N$ is almost surely finite for any integer $N$.

\vspace{0.4 cm}
The concept of metastability has now a long history in statistical physics as well as in probability theory, and because it has undergone a continuous conceptual evolution we may give some background and explain why we call our main result a result of metastability. In a seminal paper (see \cite{lebowitz}) O. Penrose and J. L. Lebowitz proposed the following characterization for a metastable thermodynamic state :

\begin{enumerate}
\item only one thermodynamic phase is present,
\item a system that starts in this state is likely to take a long time to get out,
\item once the system has gotten out, it is unlikely to return.
\end{enumerate} The first and the third points can be translated in our paradigm into the trivial fact - mentioned above - that the finite system dies out almost surely, there is therefore only one invariant measure (one phase), which is the one putting the whole mass on the state where all the neurons are quiescent. About ten years latter, M. Cassandro, A. Galves, E. Olivieri and M. E. Vares introduced in \cite{cassandro} a refinement of the second point, as they realized that a crucial characteristic of a wide variety of metastable stochastic dynamics is not only that the exit time from the metastable state is long, but that it is also in some sense unpredictable, which translate mathematically by saying that it is exponentially distributed (asymptotically). The reason for this is that the exponential distribution is characterized by the memory-less property: knowing that the system survived up to time $t$ gives you no information about what should happen next. From that point this property has been studied in a wide variety of stochastic dynamics. See for example \cite{schonmann}, \cite{mountford} or \cite{fernandez} (as a non-exhaustive list of references). Since then alternative approaches has also been developed and we refer to \cite{vares} for a complete review.

\vspace{0.4 cm}

We would also like to stress that the specific instantiation of the model we propose ourself to study is quite schematic compared to an actual neural network, especially because of our choice for the activation function $\phi$ (hard-threshold) and for the structure of the network (nearest-neighbours interaction). One of the main reason for this choice is that our purpose is to build rigorous mathematical proofs, and it should be understood that one implicit idea behind this work is that it is reasonable to expect that our main result would actually still hold for more realistic instantiations (with a sigmoid function for $\phi$, and a random graph for the network for example), even if the mathematical proof would probably be - at least for the author - out of reach in these cases.

\vspace{0.4 cm}

An important part of the present work consists in extending the result of phase transition obtained in \cite{ferrari} by proving that it holds as well for the semi-infinite version of the model (that is the process with $I=\Z^+$, the set of non-negative integers, or $I=\Z^-$ the set of non-positive integers, instead of $I=\Z$), and by deriving the various consequences that this same result has on the asymptotic distribution of the process.

\vspace{0.4 cm}

The paper is organized as follows. In Section \ref{definition} we introduce the notations and we give a proper definition of the auxiliary process. In Section \ref{dual}, we introduce the \textit{dual process}. In Section \ref{asymptotic} we obtain various results on the asymptotic distributions of the infinite process and of the semi-infinite process (by proving that the phase transition holds in the semi-infinite case as well). Finally, after establishing in Section \ref{lemmas} various relations between the infinite, semi-infinite and finite process via coupling techniques, we prove our main theorem in Section \ref{proof}.

\vspace{0.4 cm}

\section{Definition of the process}

\label{definition}

\vspace{0.4 cm}

\subsection{Formal definition via infinitesimal generator}

The stochastic process we consider is a continuous time Markov process taking values in $\{0,1\}^{\Z}$ and denoted $(\xi(t))_{t \geq 0}$. A configuration of the process is a doubly infinite sequence of $0$ and $1$ indicating in which state each neuron in the lattice is. For any $\eta \in \{0,1\}^{\Z}$, we will denote by $(\xi^\eta(t))_{t \geq 0}$ the process with initial configuration $\xi^\eta(0) = \eta$.

\vspace{0.4 cm}

Our process has the following generator:
\begin{equation} \label{generator} \L f (\eta) = \gamma \sum_{i \in \Z} \Big(f(\pi^{\dagger}_i(\eta)) - f(\eta)\Big) + \sum_{i \in \Z} \eta_i\Big(f(\pi_i(\eta)) - f(\eta)\Big),
\end{equation} where $f : \{0,1\}^{\Z} \mapsto \R$ is a cylinder function, $\gamma$ is a non-negative real number, and the $\pi^{\dagger}_i$'s and $\pi_i$'s are maps from $\{0,1\}^{\Z}$ to $\{0,1\}^{\Z}$ defined for any $i \in \Z$ as follows:

\[
    {\Big(\pi^{\dagger}_i(\eta)}\Big)_j= 
\begin{cases}
    0& \text{if } j = i,\\
    \eta_j & \text{otherwise},
\end{cases}
\] and

\[
    {\Big(\pi_i(\eta)}\Big)_j= 
\begin{cases}
    0 & \text{if } j = i,\\
    \max (\eta_i, \eta_j)  & \text{if } j \in \{i-1, i+1\},\\
    \eta_j & \text{otherwise}.
\end{cases}
\]

\vspace{0.4 cm}

It should be clear here that the $\pi^{\dagger}_i$'s correspond to the leakage effect mentioned in the informal description of the previous section, and that the $\pi_i$'s correspond to the spikes.

\vspace{0.4 cm}

For any $\eta \in \{0,1\}^\Z$ we define the \textbf{extinction time} of the processes $(\xi^\eta_t)_{t \geq 0}$ $$ \tau^\eta = \inf \{t \geq 0 : \xi^\eta(t)_i = 0 \text{ for any } i \in \Z\},$$ with the convention that $\inf \emptyset = + \infty$.
\vspace{0.4 cm}

In what follows we will often use the notation $\eta \equiv 1$ to denote the "all one" configuration and the notation $\eta \equiv 0$ to denote the "all zero" configuration. Moreover we adopt the convention of writing simply $\xi(t)$ for $\xi^\eta(t)$ and $\tau$ for  $\tau^\eta$ when the initial configuration $\eta$ is the "all one" configuration. As an abuse of notation we also write $\xi^i(t)$ for $\xi^{\{i\}}(t)$.

\vspace{0.4 cm}

\subsection{The graphical construction}

Inspired by the graphical construction introduced by Harris in \cite{harris78} we consider an alternative construction of our process.

\vspace{0.4 cm}

For any $i \in \N$ let $(N_i(t))_{t \geq 0}$ and $(N^{\dagger}_i(t))_{t \geq 0}$ be the two independent homogeneous Poisson processes mentioned in the introduction, with intensity $1$ and $\gamma$ respectively, and let $(T_{i,n})_{n \geq 0}$ and $(T^\dagger_{i,n})_{n \geq 0}$ be their respective jump times. We also impose that the collection of Poisson processes we get are defined on the same probability space $(\Omega, \F, \P)$ and are mutually independent.

\vspace{0.4 cm}

Moreover we consider the time-space diagram $\Z \times \R_+$, and for any realization of the Poisson processes, we do the following:

\begin{itemize}
\item for all $i \in \Z$ and $n \in \N$ put a "$\delta$" mark at the point $(i,T^{\dagger}_{i,n})$,
\item for all $i \in \Z$ and $n \in \N$ put an arrow pointing from $(i,T_{i,n})$ to $(i+1,T_{i,n})$ and another pointing from $(i,T_{i,n})$ to $(i-1,T_{i,n})$.
\end{itemize}

\vspace{0.4 cm}

That way we obtain a random graph $\G$ which consists of the time-space diagram $\Z \times \R$ augmented by the set of "$\delta$" marks and horizontal arrows we just described, and which is constructed on the underlying probability space $(\Omega, \F, \P)$.

\vspace{0.4 cm}

We call a \textbf{time segment} any subset of $\Z \times \R$ of the form $\{ (i,s), t \leq s \leq t'\}$, for some $i \in \Z$ and some $t < t'$. Moreover, for some $i,j \in \Z$ and $t < t'$, and for any realization of the graph $\G$, we say that there is a path from $(i,t)$ to $(j,t')$ in $\G$ if there is a connected chain of time segment and arrows leading from $(i,t)$ to $(j,t')$. We say that it is a \textbf{valid path} if it satisfies the following constraints:

\begin{itemize}
\item it never cross a "$\delta$" mark,
\item when moving upward, we never cross the rear side of an arrow.
\end{itemize}

\vspace{0.4 cm}

We write $(i,t)  \longrightarrow (j,t')$ when there is a valid path from $(i,t)$ to $(j,t')$ in $\G$. With this construction we can easily give the following characterization of our stochastic process. For any $A \in \mathcal{P}(\Z)$, and for any $t \geq 0$ :

$$\xi^A(t) = \{j \in \Z :  (i,0) \longrightarrow (j,t) \text{ for some }  i \in A\}.$$

\vspace{0.4 cm}

Notice that we moved from a process with state space $\{0,1\}^{\Z}$ to a process with state space $\mathcal{P}(\Z)$, the set of all subsets of $\Z$. It's of course only a different way to write the same thing, as any element $\eta$ of $\{0,1\}^{\Z}$ can be bijectively mapped to an element $A$ of $\mathcal{P}(\Z)$ - via the obvious relation $A = \{i \in \Z \text{ such that } \eta_i = 1\}$ - so that we can indifferently use both ways. By convention we will use $\eta, \xi \ldots$ for elements of $\{0,1\}^{\Z}$ and $A, B \ldots$ for elements of $\mathcal{P}(\Z)$. What we mean should be clear from the context.

\begin{figure}[!htb]
        \center{\includegraphics[width=13cm]
        {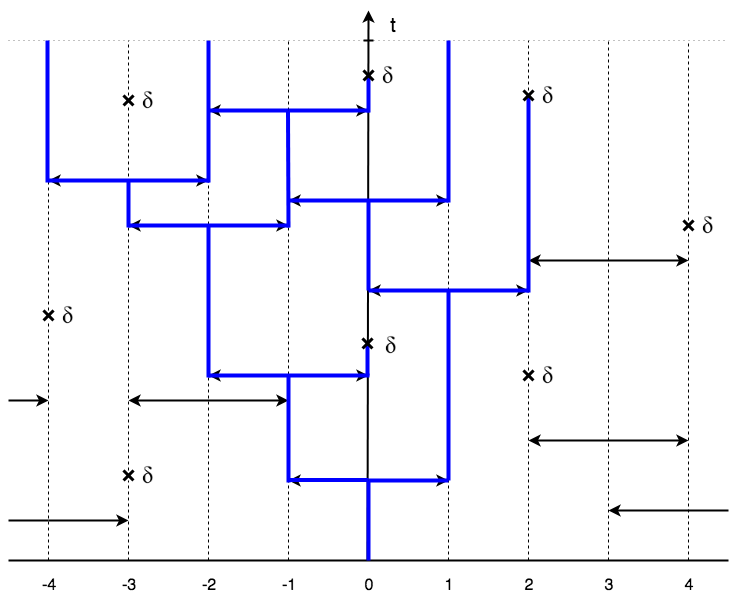}}
        \caption{\label{graph_path} In blue all the possible valid paths starting from $(0,0)$ for some realization of the graph $\G$. Here the configuration of the process at time $t$ when the initial configuration is the singleton $\{0\}$ is the set $\{-4,-2,1\}$.}
\end{figure}

\vspace{0.4 cm}

The reason for introducing this graphical construction is that it proved itself to be a powerful tool in the field of interacting particle systems.

\subsection{The finite and semi-infinite processes}

\label{finsemifin}

In order to state and prove the metastability result that we are interested in we need to introduce restricted versions of the infinite system of spiking neurons.

\vspace{0.4 cm}

For any $N \in \Z$, the \textbf{right semi-infinite process}, which we denote $(\xi_{[N,+\infty]}(t))_{t \geq 0}$, is the process taking values in $\mathcal{P}(\Z \cap [N,+\infty])$, defined as the process $(\xi(t))_{t \geq 0}$, with the random graph $\G$, but using only the $\delta$'s and arrows from the sub diagram $\{N,N+1,\ldots\} \times \R_+$. The \textbf{left semi-infinite process} $(\xi_{[-\infty,N]}(t))_{t \geq 0}$ is defined in a similar way on $\mathcal{P}(\Z \cap [-\infty,N])$, considering only $\delta$'s and arrows on the sub diagram $\{\ldots, N-1, N\} \times \R_+$. Analogously we also define the \textbf{finite process} $(\xi_N(t))_{t \geq 0}$ on $\mathcal{P}(\Z \cap [-N,N])$, considering only $\delta$'s and arrows on the sub diagram $\{-N, \ldots, N\} \times \R_+$.

\vspace{0.4 cm}

Notice that all these processes are constructed on the same probability space $(\Omega, \F, \P)$, and that they satisfy the following monotonicity relationships:

$$ \forall A \in \mathcal{P} (\Z \cap [N,+\infty]), \text{ } \forall t \geq 0, \text{ } \xi^A_{[N,+\infty]}(t) \subset \xi^A(t),$$

$$ \forall A \in \mathcal{P} (\Z \cap [-\infty,N]), \text{ } \forall t \geq 0, \text{ } \xi^A_{[-\infty,N]}(t) \subset \xi^A(t),$$

$$ \forall A \in \mathcal{P} (\Z \cap [-N,N]), \text{ } \forall t \geq 0, \text{ } \xi^A_N(t) \subset \xi^A(t).$$

\vspace{0.4 cm}

The same way we defined an extinction time for the infinite process, for any $N \in \Z$ and any initial configuration $A \in \mathcal{P} (\Z \cap [-N,N])$ we define the extinction time for the finite process, denoted $\tau^A_N$. We define as well an extinction time for the right semi-infinite process (resp. left semi-infinite process), denoted $\tau^A_{[N,+\infty[}$ (resp. $\tau^A_{]-\infty,N]}$), for any $A\in \mathcal{P} (\Z \cap [N,+ \infty])$ (resp. $A\in \mathcal{P} (\Z \cap [- \infty,N])$). We adopt the same conventions as for the infinite process regarding the notation of the initial configuration.

\vspace{0.4 cm}

Now that the objects we are interested in are well-defined and their notation clear, we state the theorem we are aimed to prove below.

\begin{theorem} \label{mainth}
There exists $\gamma'_c$ such that if $0 < \gamma < \gamma'_c$, then we have the following convergence

$$\frac{\tau_N}{\E (\tau_N)} \overset{\mathcal{D}}{\underset{N \rightarrow \infty}{\longrightarrow}} \mathcal{E} (1).$$
\end{theorem}

\vspace{0.4 cm}

\section{The dual process}

\label{dual}

\vspace{0.4 cm}

\subsection{Formal definition of the dual process}
The fact that this process is additive has the nice consequence that it has a dual. We will not explain here the details of the general definition of a dual process, and we refer to \cite{harris76} and \cite{bertein} for any reader interested in the general theory of duality.

\vspace{0.4 cm}

The \textbf{dual process} of our system of spiking neurons is defined on the state space $\mathcal{P}_f(\Z)$ of finite subset of $\Z$ and has the following generator (see \cite{ferrari}):

\begin{equation} \label{generator_dual}
\tilde{\L} g (F) = \gamma \sum_{i \in F} \Big(g(\tilde{\pi}^{\dagger}_i(F)) - g(F)\Big) + \sum_{i \in F} \eta_i\Big(g(\tilde{\pi}_i(F)) - g(F)\Big),
\end{equation} where $g$ is a cylindrical function and $F \in \mathcal{P}_f(\Z)$. The $\tilde{\pi}^{\dagger}_i$'s and $\tilde{\pi}_i$'s are called the dual maps and are defined as follows:

\[
    \tilde{\pi}^{\dagger}_i(F) = F \setminus \{i\}
\] for all $F \in \mathcal{P}_f(\Z)$ and $i \in \Z$, and

\[
    \tilde{\pi}_i(\{j\})= 
\begin{cases}
    \emptyset & \text{if } j = i,\\
    \{i, j\}  & \text{if } j \in \{i-1, i+1\},\\
    \{j\}& \text{otherwise},
\end{cases}
\] for all $i \in \Z$ and $j \in \Z$, the map for bigger sets $F \in \mathcal{P}_f(\Z)$ being given for all $i \in \Z$ by

$$ \tilde{\pi}_i(F) = \bigcup_{j \in F} \tilde{\pi}_i(\{j\}).$$

\vspace{0.4 cm}

We adopt the notation used in \cite{ferrari} and, for any $A \in \mathcal{P}_f(\Z)$, we write $(C^A(t))_{t \geq 0}$ for the process with generator (\ref{generator_dual}) and with initial configuration $C^A(0) = A$. Now, the interesting thing about duality is that the process and its dual are connected via the \textbf{duality property} (Theorem 2 in \cite{ferrari}), that we state immediately below.

\vspace{0.4 cm}

\begin{theorem} \label{duality}
For any $B \in \mathcal{P}_f(\Z)$, $A \in \mathcal{P}(\Z)$, and $t \geq 0$ we have

$$ \P \Big( \xi^A(t) \cap B \neq \emptyset \Big) = \P \Big( C^B(t) \cap A \neq \emptyset \Big).$$
\end{theorem}

\vspace{0.4 cm}

Intuitively this last theorem means that the dual process essentially behave like the initial process with the time reversed. In the following sections the duality property will often be used in the following form:

$$ \P \Big( \xi (t) \cap \{0\} \neq \emptyset \Big) = \P \Big( C^0 (t) \neq \emptyset \Big). $$

\vspace{0.4 cm}

\subsection{Graphical construction of the dual process}

\vspace{0.4 cm}

It is also possible to build a graphical construction for the dual process. Again, for any $i \in \N$ and $n \in \N$, let's consider two independent homogeneous Poisson processes $(\tilde{N}_i(t))_{t \geq 0}$ and $(\tilde{N}^\dagger_i(t))_{t \geq 0}$ with intensity $1$ and $\gamma$ respectively, and let $(\tilde{T}_{i,n})_{n \geq 0}$ and $(\tilde{T}^\dagger_{i,n})_{n \geq 0}$ be their respective jump times. As previously all the Poisson processes are assumed to be mutually independent.

\vspace{0.4 cm}

We consider the time-space diagram $\Z \times \R_+$, and for any realization of the Poisson processes, we do the following:

\begin{itemize}
\item for all $i \in \Z$ and $n \in \N$ put a "$\delta$" mark at the point $(i,\tilde{T}^{\dagger}_{i,n})$,
\item for all $i \in \Z$ and $n \in \N$ put an arrow pointing from $(i+1,\tilde{T}_{i,n})$ to $(i,\tilde{T}_{i,n})$ and another pointing from $(i-1,\tilde{T}_{i,n})$ to $(i,\tilde{T}_{i,n})$.
\end{itemize}

\vspace{0.4 cm}

As previously we get a random graph that we denote $\tilde{\G}$. Now we say that a path in $\tilde{\G}$ is a \textbf{dual-valid path} if it satisfies the following constraints:

\begin{itemize}
\item it never cross a "$\delta$" mark,
\item when moving upward, we never cross the tip of an arrow.
\end{itemize}

\vspace{0.4 cm}

We write $(i,t) \overset{\text{dual}}{\longrightarrow} (j,t')$ when there is a dual-valid path from $(i,t)$ to $(j,t')$ in $\G'$. Then, for any $A \in \mathcal{P}(\Z)$ and for any $t \geq 0$, we can write :

$$C^A(t) = \{j \in \Z :  (i,0) \overset{\text{dual}}{\longrightarrow} (j,t) \text{ for some }  i \in A\}.$$

\begin{figure}[!htb]
        \center{\includegraphics[width=13cm]
        {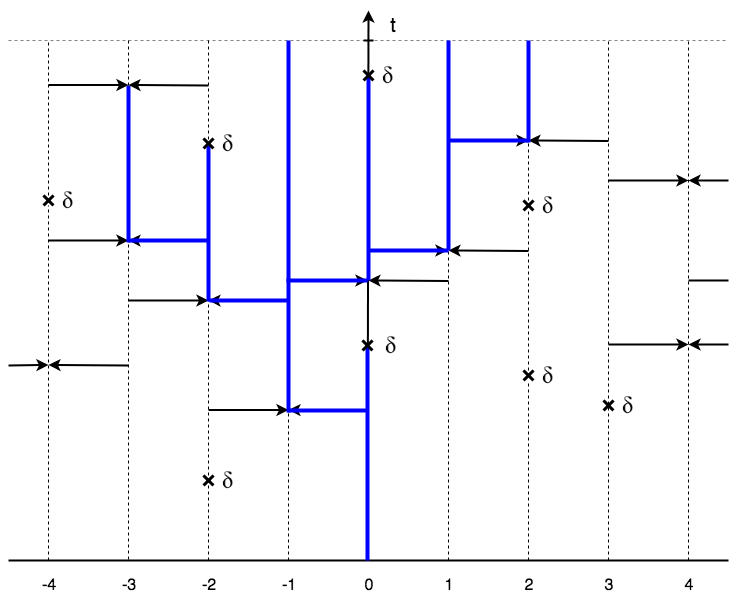}}
        \caption{\label{graphical_path_dual} In blue all the possible dual-valid paths starting from $(0,0)$ for some realization of the graph $\G'$. Here the configuration of the dual process at time $t$ when the initial configuration is the singleton $\{0\}$ is the set $\{-1,1,2\}$.}
\end{figure}

\vspace{0.5 cm}

As well as for the original process we define, for any $A \in \mathcal{P}(\Z)$, the extinction time for the dual-process:

$$ \tilde{\tau}^A = \inf \{t \geq 0 : C^A(t) = \emptyset\}.$$

The same way we defined the finite and semi-infinite processes in section \ref{finsemifin}, we define the finite and semi-infinite dual processes, using the random graph $\tilde{\G}$, and as previously we denote them $(C_{[N,+\infty[}(t))_{t \geq 0}$, $(C_{]-\infty,N]}(t))_{t \geq 0}$ and $(C_N(t))_{t \geq 0}$. We define their extinction times as well, denoted $\tilde{\tau}_N^A$, $\tilde{\tau}_{[N,+\infty[}^A$ and $\tilde{\tau}_{]-\infty,N]}^A$ for any initial configuration $A \in \mathcal{P}_f (\Z)$. We adopt the usual conventions regarding the notation for the initial configuration.

\section{Asymptotic behavior}

\label{asymptotic}

\vspace{0.4 cm}

\subsection{Asymptotic behavior of the infinite processes}

\label{sectionasinfproc}

\vspace{0.4 cm}

The dual process itself presents some kind of phase transition, as stated in the following theorem (which is Theorem 3 in \cite{ferrari}).

\begin{theorem}\label{thm:phasetrans}
There exists $0 < \gamma_c < +\infty$ such that for all $i \in \Z$ we have:

$$ \P \left(\tilde{\tau}^i = +\infty \right) > 0, \text{ if } \gamma < \gamma_c, $$

and

$$ \P \left(\tilde{\tau}^i = +\infty \right) = 0, \text{ if } \gamma > \gamma_c. $$
\end{theorem}

\vspace{0.4 cm}

A central problem we need to address is to determine the invariant measures of each of the different processes we introduced, in the sub-critical regime. Most of this analysis is done by combining Theorem \ref{thm:phasetrans} and the duality property.

\vspace{0.4 cm}

Before going any further we begin by giving a topological structure to the state space, which will be needed in what follows. We equip $\{0,1\}$ with the discrete topology so that $\{0,1\}^\Z$ can then be equipped with the corresponding product topology. That way $\{0,1\}^\Z$ is compact by Tychonoff's theorem and metrizable as any distance of the form $d(x,y) = \sum_{i\in\Z} a_i \mathbbm{1}_{\{x(i) = y(i)\}}$ generates the product topology (where $(a_i)_{i \in \Z}$ is any sequence satisfying $\sum_{i\in\Z} a_i < \infty$). The topological space $\{0,1\}^\Z$ is then associated with the corresponding Borel $\sigma$-algebra.

\vspace{0.4 cm}

Now let us introduce some order relations on $\{0,1\}^\Z$ and on the set of probability measures on $\{0,1\}^\Z$ that will be needed in what follows. Given two configurations $\eta^1$ and $\eta^2$, we will say that $\eta^1 \leq \eta^2$ if for any $i \in \Z$ we have $\eta^1_i \leq \eta^2_i$. Now, for any continuous function $f$ on $\{0,1\}^\Z$, we say that $f$ is increasing if $f(\eta^1) \leq f(\eta^2)$ whenever $\eta^1 \leq \eta^2$. We say that $f$ is decreasing if $-f$ is increasing. Finally, given two probability measures $\nu_1$ and $\nu_2$ on $\{0,1\}^\Z$, we say that $\nu_1 \leq \nu_2$ whenever the following inequality $$\int f d \nu_1 \leq \int f d \nu_2$$ holds for any continuous and increasing function $f$. One of the reasons behind this definition is that it is a well-known fact that the set of continuous and increasing function on $\{0,1\}^\Z$ is \textbf{distribution determining}, which means that for any probability measures $\nu_1$ and $\nu_2$ on $\{0,1\}^\Z$, if the following equality $$\int f d \nu_1 = \int f d \nu_2$$ holds for any continuous and increasing function $f$ then we have $\nu_1 = \nu_2$ (see the annex). In particular this implies that if $\nu_1 \leq \nu_2$ and $\nu_1 \geq \nu_2$ then $\nu_1 = \nu_2$.

\vspace{0.4 cm}

The next result - and therefore most of the results that follow - could be proved using a very general tool called the "basic coupling" (see Theorem 2.4 in Chapter 2 of \cite{ips} or section 7 of \cite{durrettintro}), nonetheless we give a somewhat more elementary proof based on the graphical construction in order to make this paper as self-contained as possible.

\begin{proposition}\label{prop:monotone}
Let $A \subset B \subset \Z$, then for any $t \geq 0$ we have

$$\xi^A(t) \subset \xi^B(t).$$

Moreover, for any $0 \leq s < t$ we have the following
$$\P \left( \xi(s) \in \bigcdot \text{ }\right) \geq \P \left( \xi(t) \in \bigcdot \text{ }\right).$$
\end{proposition}

\begin{proof}
For the first part of the proposition it suffices to notice that if $(i,0) \longrightarrow (j,t)$ for some $i$ in $A$ and some $t \geq 0$ then $(j,t) \in \xi^B(t)$ as $i$ belongs to $B$ as-well. For the second part of the proposition fix $0 \leq s < t$. Then let $\eta$ be a random variable taking value in $\{0,1\}^\Z$ and having the same distribution as $\xi(t-s)$. By the first part we have $\xi(s) \geq \xi^\eta (s)$, so that for any continuous and increasing function $f$ we have $f(\xi(s)) \geq f(\xi^\eta (s))$, and taking the expectation we get $\E \left( f(\xi(s)) \right) \geq \E \left( f(\xi^\eta (s)) \right)$. But by construction $\xi^\eta (s)$ has the same distribution as $\xi(t)$ so that we end up with  $$ \E \Big( f(\xi(s)) \Big) \geq \E \Big( f(\xi (t)) \Big),$$ which is the same as $$ \int f d\P(\xi(s) \in \bigcdot \text{ }) \geq \int f d\P(\xi(t) \in \bigcdot \text{ }).$$
\end{proof}

\begin{remark}
To avoid confusion we call the first property \textbf{set monotonicity}, and the second one \textbf{stochastic monotonicity}.
\end{remark}

\vspace{0.4 cm}

From this last proposition we get the following.

\begin{corollary} \label{prop:asymptomeasure}
For any $\gamma > 0$ there exists a probability measure $\mu_\gamma$ that is invariant for $(\xi(t))_{t \geq 0}$ and that is such that $$\P \Big( \xi(t) \in \bigcdot \text{ } \Big) \underset{t \rightarrow \infty}{\longrightarrow} \mu_\gamma.$$ The Dirac measure on the "all zero" configuration, denoted $\delta_\emptyset$, is also invariant and we have $$\P \Big( \xi^\emptyset(t) \in \bigcdot \text{ } \Big) \underset{t \rightarrow \infty}{\longrightarrow} \delta_\emptyset.$$ Moreover if $\nu$ is any other invariant measure then $\delta_\emptyset \leq \nu \leq \mu_\gamma$.
\end{corollary}

\begin{proof}
The fact that the process starting from the "all zero" configuration converges weakly to $\delta_\emptyset$ is of course entirely trivial as we actually have $\xi^\emptyset(t) = \emptyset$ for any $t \geq 0$. For the convergence of the process starting from the "all one" configuration take any continuous and increasing function $f$ and remember that we showed in the proof of the previous proposition that $t \mapsto \E \big(f(\xi(t))\big)$ is a decreasing function. It follows that  $\E \big(f(\xi(t))\big)$ converges to some finite constant when $t$ goes to infinity (remember that $\{0,1\}^\Z$ is compact). The set of continuous and increasing functions being distribution determining the result follows.

For the last statement of the proposition suppose we have some invariant measure $\nu$ and write $(\xi^\nu(t))_{t \geq 0}$ for the process with initial configuration chosen randomly with respect to distribution $\nu$. For any $t \geq 0$ by set monotonicity we have $\xi^\emptyset(t) \leq \xi^\nu(t) \leq \xi(t)$, so for any continuous and increasing function $f$ we have 

$$\E \Big(f(\xi^\emptyset(t))\Big) \leq \E \Big(f(\xi^\nu(t))\Big) \leq \E \Big(f(\xi(t))\Big),$$ or equivalently $$ \P \Big(\xi^\emptyset(t) \in \bigcdot \text{ }\Big) \leq \P\Big(\xi^\nu(t) \in \bigcdot \text{ }\Big) \leq \P\Big(\xi(t) \in \bigcdot \text{ }\Big).$$

The inequality $\delta_\emptyset \leq \nu \leq \mu_\gamma$ then follows from the convergence results proven above and from the fact that $\P\big(\xi^\nu(t) \in \bigcdot \text{ }\big) = \nu$ for any $t \geq 0$.
\end{proof}

\vspace{0.4 cm}

The asymptotic distribution of the process starting from the "all one" configuration will be referred as the \textbf{upper-invariant measure} and the asymptotic distribution of the process starting from the "all zero" configuration will be referred as the \textbf{lower-invariant measure}. The dual process has an upper-invariant measure too, which we denote $\tilde{\mu}_\gamma$, and his lower-invariant measure is also the Dirac $\delta_\emptyset$ (to see this it suffices to verify that all the arguments used above work for the dual process as well). Moreover the inequality $\delta_\emptyset \leq \nu \leq \tilde{\mu}_\gamma$ remains true if $\nu$ is an invariant measure for the dual process.

\vspace{0.4 cm}

The fact that the upper-invariant and lower-invariant measures are a lower and upper bound respectively for any invariant measure has the following consequence.

\vspace{0.4 cm}

Let define the \textbf{density} of the process $(\xi(t))_{t \geq 0}$:

$$\rho_\gamma = \mu_\gamma \left( \{\eta : \eta_0 = 1 \} \right).$$

\vspace{0.4 cm}

As a consequence of Theorem \ref{thm:phasetrans} we have the following result for the density.

\begin{proposition} \label{prop:density}
When $\gamma < \gamma_c$ we have $\rho_\gamma > 0$, and therefore $\mu_\gamma \neq \delta_\emptyset$.
\end{proposition}

\begin{proof}
This follows from duality (Theorem \ref{duality}) and can be derived as follows

\begin{align*}
    \rho_\gamma &= \lim_{t \rightarrow \infty} \P \Big( \xi (t) \cap \{0\} \neq \emptyset \Big)\\
    &= \lim_{t \rightarrow \infty} \P \Big( C^0 (t) \neq \emptyset \Big)\\
    &= \lim_{t \rightarrow \infty} \P \left( \tilde{\tau}^0 > t\right)\\
    &= \P \left( \tilde{\tau}^0 = \infty\right),
\end{align*} and $\P \left( \tilde{\tau}^0 = \infty\right) > 0$ when $\gamma > \gamma_c$.
\end{proof}

\vspace{0.4 cm}

In order to prove the metastability result we are only interested in the sub-critical regime so that from now on we will just assume that $\gamma < \gamma_c$ and omit the dependence in $\gamma$ in the notation, writing simply $\mu$ for $\mu_\gamma$, $\tilde{\mu}$ for $\tilde{\mu}_\gamma$, and $\rho$ for $\rho_\gamma$. We have the following result, which is the equivalent for the dual process of the second part of Proposition \ref{prop:density}.

\begin{proposition}
In the sub-critical regime $\tilde{\mu} \neq \delta_\emptyset$.
\end{proposition}

\begin{proof}
\begin{align*}
    \tilde{\mu} (\eta \equiv 0) &= \lim_{t \rightarrow \infty} \P \left( C (t) = \emptyset \right)\\
    &= 1 - \lim_{t \rightarrow \infty} \P \left( C (t) \neq \emptyset \right)\\
    &\leq  1 - \lim_{t \rightarrow \infty} \P \left( C^0 (t) \neq \emptyset \right)\\
    &= 1 - \P \left( \tilde{\tau}^0 = \infty\right),
\end{align*} so that $\tilde{\mu} (\eta \equiv 0) < 1$.
\end{proof}

\vspace{0.4 cm}

To state the lemma below let introduce the following notation: $\I$ (resp. $\tilde{\I}$) will denote the set of invariant measures of the process $(\xi(t))_{t \geq 0}$ (resp. $(C(t))_{t \geq 0}$). We know, by classical theory of Markov processes (see for example proposition 1.8 of chapter 1 in \cite{ips}) that $\I$ and $\tilde{\I}$ are convex sets, so that we can define $\I_e$ and $\tilde{\I}_e$ the set of extreme points of $\I$ and $\tilde{\I}$ respectively, and we know that $\I$ and $\tilde{\I}$ are the convex-hull of $\I_e$ and $\tilde{\I}_e$ respectively (as a consequence of the Krein-Milman theorem).

\begin{lemma} \label{lemma:extremal}
We have $\{\delta_\emptyset,\mu\} \subset \I_e$ and $\{\delta_\emptyset,\tilde{\mu}\} \subset \tilde{\I}_e$. In words, the upper-invariant and lower-invariant measures are extremal.
\end{lemma}

\begin{proof}
We prove the statement for $\mu$. Suppose that there exists $\nu_1$ and $\nu_2$ in $\I$ such that $\mu = p \nu_1 + (1-p) \nu_2$ for some $0<p<1$. Then Proposition \ref{prop:asymptomeasure} gives us that $\mu_1 \leq \mu$ and $\mu_2 \leq \mu$ so that for any continuous and monotone function $f$ we have

$$\int f d\nu_1 \leq \int f d\mu \text{ and } \int f d\nu_2 \leq \int f d\mu.$$

But we also have $$\int f d\mu = p \int f d\nu_1 + (1-p) \int f d\nu_2,$$ from what if follows that $$\int f d\mu = \int f d\nu_1 = \int f d\nu_2.$$

\end{proof}

\vspace{0.4 cm}

Using this last lemma it can be shown that not only $\mu$ (as well as $\tilde{\mu}$) is different from $\delta_\emptyset$ in the sub-critical regime, but it also puts no mass on $\eta \equiv 0$. This is the subject of the following proposition.

\begin{proposition}\label{prop:mumasszero}
In the sub-critical regime we have $\mu \left( \eta \equiv 0\right) = 0$ and $\tilde{\mu} \left( \eta \equiv 0\right) = 0$.
\end{proposition}

\begin{proof}
We prove the statement for $\mu$. Regardless of the value of  $\mu \left( \eta \equiv 0\right)$ we can always find some $p \in [0,1]$ and some probability measure $\nu$ satisfying $\nu(\eta \equiv 0) = 0$ such that

$$\mu = p \delta_\emptyset + (1-p) \nu.$$

We know that $\mu \neq \delta_\emptyset$ so that $p$ has to be different from $1$, thus $\nu$ need to be invariant as well. But if $p$ were different from $0$ then $\mu$ would be a (non-trivial) convex combination of invariant measures, which would be a contradiction with Lemma \ref{lemma:extremal}. We conclude that $\mu = \nu$.
\end{proof}

\vspace{0.4 cm}

Finally, one important result we will need in order to prove metastability is the spatial ergodicity of the measure $\mu$. It is stated in the following theorem.

\begin{theorem} \label{ergodicity}
The measure $\mu$ is spatially ergodic in the sense that a sequence of random variable $(X_k)_{k\in \Z}$ taking value in $\{0,1\}$ and such that $X_k$ is distributed like $\mu \big( \{\eta: \eta_k = \bigcdot \text{ }\}\big)$ would satisfy the following

$$ \frac{1}{n+1} \sum_{k=0}^n X_k \overset{a.s.}{\underset{n \rightarrow \infty}{\longrightarrow}} \rho.$$
\end{theorem}

\begin{proof}
Using a similar coupling as in the proof of Proposition \ref{prop:monotone} we can construct an infinite sequence of random variables in $\{0,1\}^\Z$, denoted $\left(\eta^k \right)_{k \in \N \cup \{\infty\}}$, satisfying $\eta^0 \geq \eta^1 \geq \eta^2 \geq \ldots \geq \eta^\infty$ and such that $\eta^0$ is equal to $\eta \equiv 1$, $\eta^k$ has the same distribution as $\xi(k)$ for any $k \geq 0$, and $\eta^\infty$ is distributed according to $\mu$. Let $\theta$ be the shift operator on $\{0,1\}^\Z$, i.e. the operator defined for any $\eta \in \{0,1\}^\Z$ and $x \in \Z$ by $(\theta \eta )(x) = \eta (x + 1)$. For $k \geq 0$ the composition of order $k$ of $\theta$ with itself will be denoted $\theta^k$. Moreover let $\mathbbm{1}_0$ be the function defined for any $\eta \in \{0,1\}^\Z$ by $\mathbbm{1}_0 (\eta) = \mathbbm{1}_{\{\eta_0 = 1\}}$. Then for any $m\geq0$ and $n\geq0$ we have 

\begin{equation} \label{ergosums}
\frac{1}{n+1} \sum_{k=0}^n \mathbbm{1}_0 ( \theta^k (\eta^m)) \geq  \frac{1}{n+1} \sum_{k=0}^n \mathbbm{1}_0 ( \theta^k (\eta^\infty)).
\end{equation}

\vspace{0.4 cm}

For any $t \geq 0$ we know that $(\xi_k (t))_{k \in \Z}$ is an ergodic stationary sequence (this is true for any system with finite range interaction, see \cite{holley} page 1967). Therefore, if we denote by $\mathcal{S}$ the sigma-algebra of shift invariant events with respect to $\mu$ (see Chapter 7 of \cite{tande} for precise definitions) then by Birkhoff's ergodic theorem the right-hand side of (\ref{ergosums}) converges to $\E \big(\mathbbm{1}_0(\eta^\infty) \text{ } | \text{ } \mathcal{S} \big)$ almost surely when $n$ goes to $\infty$ while the left-hand side converges to $\E \big(\mathbbm{1}_0(\eta^m)\big)$. It follows that for any $m \geq 0$ we have

$$ \E \Big(\mathbbm{1}_0(\eta^m)\Big) \geq  \E \Big(\mathbbm{1}_0(\eta^\infty) \text{ } | \text{ } \mathcal{S} \Big) \text{  a.s.}$$

And taking the limit when $m$ goes to $\infty$

$$ \E \Big(\mathbbm{1}_0(\eta^\infty)\Big) \geq  \E \Big(\mathbbm{1}_0(\eta^\infty) \text{ } | \text{ } \mathcal{S} \Big) \text{  a.s.}$$

But a real-valued random variable which is bounded by its own expectation need to be almost surely equal to it, so that in the end we have $$\frac{1}{n+1} \sum_{k=0}^n \mathbbm{1}_0 ( \theta^k (\eta^\infty)) \overset{a.s.}{\underset{n \rightarrow \infty}{\longrightarrow}} \E \Big(\mathbbm{1}_0(\eta^\infty)\Big),$$ and it suffices to point out that $ \E \big(\mathbbm{1}_0(\eta^\infty)\big) = \rho$ to end the proof.

\end{proof}

\subsection{Asymptotic behavior of the semi-infinite processes}

In order to show that the asymptotic behavior of the semi-infinite processes is essentially the same as the asymptotic behavior of the infinite process, we need to make sure that we have an equivalent of Theorem \ref{thm:phasetrans} for the semi-infinite dual process, so that the developments of section \ref{sectionasinfproc} remain valid in the semi-infinite case. This question of whether or not the phase transition remains true, and if it does, with the same critical value, is indeed not trivial. One could indeed imagine that the boundary on the left or on the right somewhat produce a different behavior. Moreover, as it will appear later, it is a crucial point for the proof of Theorem \ref{mainth} that the phase transition remains true for the semi-infinite processes. The proof uses a contour argument (see for example \cite{griffeath}) and is somewhat similar to the proof presented in \cite{ferrari} for the infinite process, nonetheless this former proof was quite elliptical and one of the goals we're pursuing here is to give a clearer argument. We observe that the phase transition for the original semi-infinite process (in the form of Theorem \ref{thm:phasetransition}) is a consequence of the Theorem \ref{theorem:semiinfphasetrans} proven below—which is only concerned with the dual process—using the same standard arguments as in section 4 of \cite{ferrari}. Nonetheless we don't state this result in this form here as it is unnecessary for our main purpose, which is metastability.

\vspace{0.4 cm}

Notice that the following Theorem is stated for the process defined on $\mathcal{P}(\Z \cap [0,+\infty[)$ but that by symmetry it obviously remains true for the process defined on $\mathcal{P}(\Z \cap [-\infty,0])$, and more generally for the processes defined on $\mathcal{P}(\Z \cap [N,+\infty[)$ or $\mathcal{P}(\Z \cap [-\infty,N])$ for any value of $N \in \Z$.

\vspace{0.4 cm}

In order to prove this phase transition, we need a preliminary result of monotonicity in $\gamma$. In the following lemma, we write $C^{i,\gamma}_{[0,\infty]}$ to make explicit the dependence in $\gamma$.

\begin{lemma}\label{lemma:lambamon}
For any $i \in  \Z \cap [0,+\infty]$ and for any $t \geq 0$ we have the following:

$$ \text{if } \gamma_1 < \gamma_2 \text{, then } \P \left( C^{i,\gamma_1}_{[0,\infty]} (t) = \emptyset\right) \leq \P \left( C^{i,\gamma_2}_{[0,\infty]} (t) = \emptyset\right).$$
\end{lemma}

\begin{proof}
The proof is exactly the same as the proof of Lemma 5 in \cite{ferrari}.
\end{proof}

\vspace{0.4 cm}

\begin{theorem} \label{theorem:semiinfphasetrans}
There exists $0 < \gamma'_c < +\infty$ such that for all $i \in \Z \cap [0,+\infty]$ we have: $$ \P \left(\tilde{\tau}^i_{[0,+\infty[} = +\infty \right) > 0, \text{ if } \gamma < \gamma'_c, $$ and $$ \P \left(\tilde{\tau}^i_{[0,+\infty[} = +\infty \right) = 0, \text{ if } \gamma > \gamma'_c. $$
\end{theorem}

\begin{proof}
Notice that by Lemma \ref{lemma:lambamon} the function $\gamma \mapsto \P \left(\lim_{t \rightarrow \infty}  C^{0,\gamma}_{[0,\infty]} (t) = \emptyset \right)$ is non-decreasing, so that in order to prove the theorem it is sufficient to find two different values of the parameter $\gamma$ such that $\P\left(\tilde{\tau}^i_{[0,+\infty[} = +\infty \right) > 0$ for any $i \in \Z$ for the first one and $\P \left(\tilde{\tau}^i_{[0,+\infty[} = +\infty \right) = 0$ for the second one.

\vspace{0.4 cm}

The second part of the Theorem is immediate by monotonicity. Indeed, for $\gamma > \gamma_c$ (where $\gamma_c$ is the critical value for the infinite process) and for any $i \in \Z$

$$ \P \left(\tilde{\tau}^i_{[0,+\infty[} = +\infty \right) \leq \P \left(\tilde{\tau}^i = +\infty \right) = 0.$$

Moreover it has been proven in \cite{ferrari} that $\gamma_c < 1$, so that we also have $\gamma'_c <1$.

\vspace{0.4 cm}

To fix ideas, we prove the first part of the Theorem for the process starting at $\{0\}$. The general result will then follow as, for any $i \in \Z^+$, we have $\P \left(\tilde{\tau}^i_{[0,+\infty[} = +\infty \right) \geq \P \left(\tilde{\tau}^0_{[0,+\infty[} = +\infty \right)$ (this can be showed by coupling, using a modification of the graph $\tilde{\G}$ where all events has been shifted spatially to the right $i$ times).

\vspace{0.4 cm}

For the first part of the theorem, we will use the fact that the event $\{\tilde{\tau}^0_{[0,+\infty[} < \infty\}$ is equivalent to the event that $C^0 \text{ } \mydef \text{ } \bigcup_{t \geq 0} C^0_{[0,+\infty[} (t)$ is a finite set. It will therefore be sufficient to prove that $\P \left( | C^0 | < \infty \right) < 1$. In order to do this we consider a realization of $| C^0 | < \infty$ and draw its contour by embedding the time-space diagram $\Z \times \R_+$ in $\R \times \R_+$ and by then defining

$$ E \text{ } \mydef \text{ } \Big\{(y,t) : |y-j| \leq \frac{1}{2} \text{ for some } j \in C^i(t), \text{ } t \geq 0\Big\}.$$

Now let fill the holes of $E$ to get the set $\tilde{E}$, and let $\Gamma$ be the boundaries of $\tilde{E}$. $\Gamma$ consists of a sequence of alternating horizontal and vertical segments and with a little thought is should be clear that there are exactly $4n$ of them (for some $n \in \Z$). Moreover we encode $\Gamma$ not as a sequence of horizontal and vertical segments but as a sequence of direction triplets $(D_1, D_2, \ldots D_{2n})$. This is done as follows: start at $(\frac{1}{2},0)$ and follow the boundary of $\Gamma$ in counterclockwise direction, at step $i$ you'll be going trough $D_i$, which is one of the seven possible triplets: $$uru, ulu, uld, drd, dru, dld, dlu,$$ where $u$, $l$, $d$, and $r$ means "up", "left", "down", and "right". The last direction of the current triplet is the first direction of the next one.

\vspace{0.4 cm}

\begin{figure}[!htb]
        \center{\includegraphics[width=9cm]
        {contour.png}}
        \caption{Example of a possible contour with $n = 7$. Some of the direction vectors are explicitly drawn.}
\end{figure}

\vspace{0.4 cm}
Now we need to bound the probability of occurrence of some of these events. For reasons that will become clear soon it is sufficient to bound the occurrence of $dld$, $dlu$ and $ulu$.

\vspace{0.4 cm}

First we find a bound for the probability of occurrence of both $dld$ and $dlu$. For a given $j \in \{1, \ldots 2n\}$, let consider what could happen to $D_j$. Consider the first horizontal segment which is immediately before the first vertical segment of $D_j$ (see figure \ref{phasetrans}), there is two possibilities:

\begin{itemize}
    \item if it is oriented to the left, let $(k,t) \in \Z \times \R_+$ be the coordinates of the point immediately to the left of its left extremity,
    \item if it is oriented to the right, let $(k,t) \in \Z \times \R_+$ be the coordinates of its midpoint.
\end{itemize}

\vspace{0.4 cm}

Notice that there is one possible case in which $k=-1$, when we are hitting the left border of our restricted space-time diagram, but in this particular case the occurrence of a $dlu$ or $dld$ is not possible, so that we can simply assume that $k \geq 0$. Now let define $$ F^{\dagger}_j = \sup \{ s \leq t : \tilde{T}^{\dagger}_{k,n} = s, \text{ } n \geq 0 \},$$ and $$ F_j = \sup \{ s \leq t : \tilde{T}_{k,n} = s, \text{ } n \geq 0 \}.$$

\vspace{0.4 cm}

It is not hard to see that (see again figure \ref{phasetrans}) $$ \{ D_j = dlu\} \cup \{ D_j = dld\} \subset \{ F^{\dagger}_j \geq F_j \} \text{ } \mydef \text{ } E_j.$$

\vspace{0.4 cm}

In words, the occurrence of $dld$ or $dlu$ implies that when you follow the line $\{k\} \times \R_+$ downward starting from $(k,t)$, then the first event you'll encounter is a $\delta$ (see figure \ref{phasetrans}).

\vspace{0.4 cm}

Finally, as it doesn't matter for the jump times of a Poisson process whether the time goes upward or downward, we have

$$ \P \left( E_j \right) = \frac{\gamma}{1+\gamma} \leq \gamma.$$

\vspace{0.4 cm}

It remains to bound the occurrence of $ulu$. using the same notation as in the precedent cases for the coordinates of the objects considered, let's define $$ G^{\dagger}_j = \sup \{ s \geq t : \tilde{T}^{\dagger}_{k,n} = s, \text{ } n \geq 0 \},$$ and $$ G_j = \sup \{ s \geq t : \tilde{T}_{k+1,n} = s, \text{ } n \geq 0 \}.$$

\vspace{0.4 cm}

Then we have

$$ \{D_j = ulu\} \subset \{G^\dagger_j \leq G_j\} \text{ } \mydef \text{ } E'_j.$$

Indeed, the occurrence of $ulu$ imply that when you follow simultaneously the two lines $\{k\}\times\R_+$ and $\{k+1\}\times\R_+$ upward, starting from $(k,t)$, you'll encounter a $\delta$ on the line $\{k\}\times\R_+$ before you encounter a spike on the line $\{k+1\}\times\R_+$ (see again figure \ref{phasetrans}). Moreover we have 

$$\P \left( E'_j \right) = \frac{\gamma}{1+\gamma} \leq \gamma.$$

\begin{figure}[!htb]
        \center{\includegraphics[width=\textwidth]
        {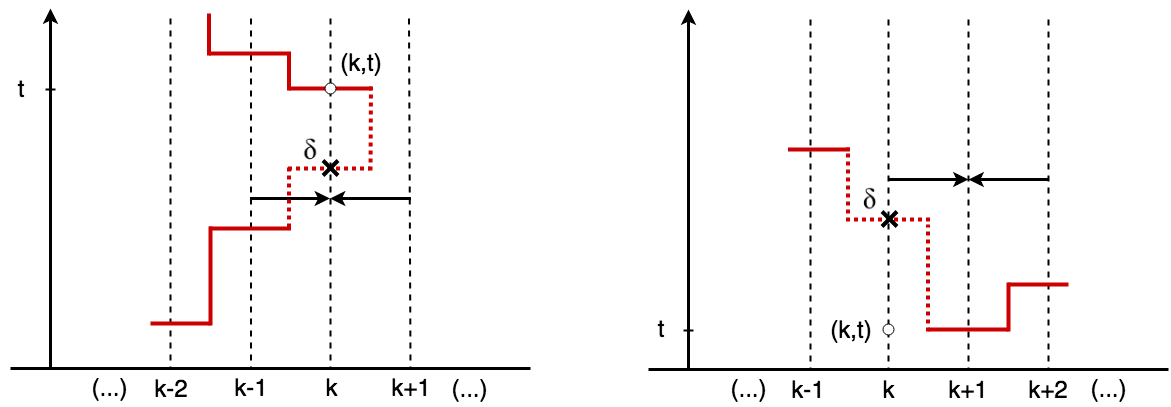}}
        \caption{\label{phasetrans} On the left side we can see $F^{\dagger}_j$ and $F_j$ that are represented on the time-space diagram. The $\delta$ corresponds to $F^{\dagger}_j$ and the double arrow corresponds to $F_j$. $G^{\dagger}_j$ and $G_j$ are represented on the right side. The red line represents $\Gamma$ and the dashed part of it represents $D_j$, which took value $dld$ on the left, and $ulu$ on the right.}
\end{figure}

\vspace{0.4 cm}

Now we would like to use the mutual independence of all the events $E_j$ and $E'_j$ in order to get a bound for our contour. Unfortunately there is some cases in which these events are actually not independent. Indeed, any event $E'_j$ (corresponding to an $ulu$) makes use of both the deaths occurring on the concerned  portion of some line $k \times \R^+$ and of the double arrows occurring on the concerned portion of the line $(k+1) \times \R^+$. However, if the preceding triplet is a $dlu$ (as in the right part of Figure \ref{phasetrans}), then a portion of the line $(k+1) \times \R^+$ that was used in order to bound the $ulu$ is also used to bound this $dlu$, so that the independence doesn't hold in this specific case. We solve this problem by avoiding contiguity using a partition of the $D_i$'s in two subsets.

\vspace{0.4 cm}

A given triplet $D_j$ crosses one single line of the form $k \times \R^+$. If $k$ is even we say that $D_j$ is of type 1 whereas if $k$ is odd we say that $D_j$ is of type 2. That way all $D_i$'s of a given type are non-contiguous (i.e. depend on different Poisson processes or at least on disjoint regions of the same Poisson processes).

\vspace{0.4 cm}

Now we write $N(dld)$ for the number of occurrences of $dld$ in a given contour, $N(dlu)$ for the number of occurrences of $dlu$ and so on. We also write $N_1(dld)$ and $N_2(dld)$ for the number of $dld$ of type 1 and type 2 respectively, $N_1(dlu)$ and $N_2(dlu)$ for $dlu$ and so on.

\vspace{0.4 cm}

For a given  contour $\Gamma$, from the discussion above and from the fact that we can just discard the occurrences of the other triplets from the intersection of events in which $\Gamma$ consists, it follows that we can bound its probability by $$\gamma^{N_1(dld) + N_1(dlu) + N_1(ulu)},$$ or indifferently by $$\gamma^{N_2(dld) + N_2(dlu) + N_2(ulu)}.$$

\vspace{0.4 cm}

Now notice that $\Gamma$ necessarily contains the same number of left and right oriented segments so that in particular, as $\Gamma$ contains $2n$ horizontal segments, it shall contain exactly $n$ segments oriented to the left. Thus we have the following equation

$$ N(dld) + N(dlu) + N(uld) + N(ulu) = n.$$

\vspace{0.4 cm}

Moreover it is not difficult to see that $$ N(uld) = N(dlu) + 1,$$ which, together with the previous equation, allows the following bound

$$ N(uld) \leq \frac{n+1}{2}.$$

\vspace{0.4 cm}

It follows that

$$ N(dld) + N(dlu) + N(ulu) \geq \frac{n-1}{2},$$

so that we either have

$$ N_1(dld) + N_1(dlu) + N_1(ulu) \geq \frac{n-1}{4},$$

or

$$ N_2(dld) + N_2(dlu) + N_2(ulu) \geq \frac{n-1}{4}.$$

\vspace{0.4 cm}

In both cases we get the following bound

$$ \P\Big( \Gamma \Big) \leq (\sqrt[4]{\gamma})^{n-1}.$$

\vspace{0.4 cm}

For $n=1$, the probability of having a contour of length $4$ is $\frac{1+\gamma}{2+\gamma}$.

\vspace{0.4 cm}

For $n=2$, the $2$ possibilities for the shape of $\Gamma$ can be bounded by $\gamma$ as both of them contain at least an $ulu$ or a $dld$, so that the probability of having a contour of length $8$ can be bounded by $2 \gamma$.

\vspace{0.4 cm}

Observe that the number of possible triplets for an arbitrary $n$—remembering that a given $\Gamma$ contains $2n$ triplets—can be bounded by $4^{2n} = 16^n$. Indeed the first segment of the first triplet is always oriented upward, and the first segment of any other triplet is determined by the last segment of the previous one, so that the number of possibilities can be roughly bounded by $4$ for each of the $D_j$'s.

\vspace{0.4 cm}

Finally, for $\gamma < \frac{1}{16^4}$, we get the following bound

\begin{align*}
 \P \left( \tau^0 < \infty \right) &\leq \frac{1 + \gamma}{2 + \gamma} + 2 \gamma + \sum_{n \geq 3} 16^n (\sqrt[4]{\gamma})^{n-1}\\
 &= \frac{1 + \gamma}{2 + \gamma} + 2 \gamma + 16^3 \sqrt{\gamma} \cdot \frac{16 \sqrt[4]{\gamma}}{1 - 16 \sqrt[4]{\gamma}}.
\end{align*}

\vspace{0.4 cm}

When $\gamma$ goes to $0$ the right-hand side of the latter inequality goes to $\frac{1}{2} < 1$, from what it follows that $\gamma'_c > 0$.

\end{proof}

\vspace{0.4 cm}

Having proven the phase transition for the semi-infinite dual process we find ourself in the same situation as we were for the infinite process at the beginning of section \ref{sectionasinfproc}, and it is easy to check that all the arguments given there remain valid in the semi-infinite case. The only results that we really are interested in are the ones concerning the sub-critical regime so we assume $\gamma < \gamma'_c$. We have stochastic monotonicity and we can define 

$$\mu_{[0,+\infty[} \text{ } \mydef \text{ } \lim_{t \rightarrow \infty} \P \left( \xi_{[0,+\infty[}(t) \in \bigcdot \text{ } \right),$$ and $$\tilde{\mu}_{[0,+\infty[} \text{ } \mydef \text{ } \lim_{t \rightarrow \infty} \P \left( C_{[0,+\infty[}(t) \in \bigcdot \text{ } \right).$$

\vspace{0.4 cm}

Moreover $\mu_{[0,+\infty[}$ and $\tilde{\mu}_{[0,+\infty[}$ are extremal invariant, and in the sub-critical regime we have that $\mu_{[0,+\infty[} \neq \delta_\emptyset$ and $\tilde{\mu}_{[0,+\infty[} \neq \delta_\emptyset$, from what it follows that we actually have 
$$\mu_{[0,+\infty[} (\eta \equiv 0) = 0 \text{ } \text{ } \text{ } \text{ } \text{and} \text{ } \text{ } \text{ } \text{ }  \tilde{\mu}_{[0,+\infty[} (\eta \equiv 0) = 0.$$

\vspace{0.4 cm}

These few facts should be remembered as they will play an important role in the proof of our main theorem.

\vspace{0.5 cm}

\section{Some technical lemmas}

\label{lemmas}

Before entering the discussion about metastability, we establish four lemmas that will be needed in the course of the proof. The first three are only consequences of the nearest-neighbours nature of the interaction, even if it might not be immediately clear from the proofs, which entirely rely on the graphical construction. The last one is simply a rigorous statement of an intuitive fact, namely that the more scattered your initial configuration is the higher is the probability to survive. The reader in a hurry might simply skip this part to come back to it later if needed.

\vspace{0.4 cm}

\begin{lemma} \label{xionlnrn}
Define $r_N(t) = \max \xi_N(t)$ and $l_N(t) = \min \xi_N(t)$. Then, for any $0 \leq t \leq \tau_N$, we have the following

$$ \xi_N(t) \cap [l_N(t),r_N(t)] = \xi(t) \cap [l_N(t),r_N(t)].$$
\end{lemma}

\begin{proof}
We only need to show that $\xi(t) \cap [l_N(t),r_N(t)] \subset \xi_N(t) \cap [l_N(t),r_N(t)]$. Let $x \in \xi(t) \cap [l_N(t),r_N(t)]$. As $x \in \xi(t)$ there exists $y \in \Z$ such that there is a valid path in the graph $\G$ from $(y,0)$ to $(x,t)$, and we denote it $P_{y \rightarrow x}$. Now we define the left and right frontiers of the finite process

$$ \partial^{\text{left}}_N (t) \text{ } \mydef \text{ } \{ (l_N(s),s) \in \Z \times \R_+, \text{ } s \in [0,t]\}$$ and $$ \partial^{\text{right}}_N (t) \text{ } \mydef \text{ } \{ (r_N(s),s) \in \Z \times \R_+, \text{ } s \in [0,t]\}.$$

\vspace{0.4 cm}

Assume first that $y \not\in \{-N, \ldots N\}$. Then the path $P_{y \rightarrow x}$ needs to cross at least one of the frontier. Let $t'$ ($t'<t$) be the last time of crossing and without loss of generality assume that this crossing is a crossing of the left frontier. Then $P^{t',t}_{y \rightarrow x}$ - the restriction of $P_{y \rightarrow x}$ to the time interval $[t',t]$ - is a valid path from $(x',t') = \partial^{\text{left}}_N (t')$ to $(x,t)$, but by definition $x' \in \xi_N(t')$ so that for some $y' \in \{-N,\ldots N\}$ there is a valid path $Q_{y' \rightarrow x'}$ from $(y',0)$ to $(x',t')$. Finally the concatenation of $Q_{y' \rightarrow x'}$ and $P^{t',t}_{y \rightarrow x}$ is a valid path from $(y',0)$ to $(x,t)$, which prove that $x \in \xi_N(t)$. Note that in the two last sentences when we used the expression "valid path" what we really meant is valid path for the finite process.

\vspace{0.4 cm}

If $y \in \{-N, \ldots N\}$, then either $P_{y \rightarrow x}$ stays inside the frontiers of the finite process and there is nothing to prove, either it crosses one of the frontiers and the argument is the same as above.
\end{proof}

\vspace{0.4 cm}

\begin{lemma} \label{lemma:rl}
Fix some $N \in \N^*$ and for some non-empty sets $B \subset \Z \cap [1,N]$ and $C \subset \Z \cap [-N,-1]$ define the stopping times $$R^B_N = \inf \{ t > 0: -N \in \xi^B_{]-\infty,N]} (t) \},$$ and $$L^C_N = \inf \{ t > 0: N \in \xi^C_{[-N,\infty[} (t) \}.$$

\vspace{0.4 cm}

If $\tau^{B \cup C}_N > \max \left(R^B_N,L^C_N \right)$, then for any $t > \max \left(R^B_N,L^C_N \right)$,

$$\xi_N (t) = \xi^{B \cup C}_N (t).$$

\end{lemma}

\begin{proof}
The proof is similar to the proof of the previous lemma and relies essentially on the fact that the interaction is between nearest neighbours. We assume $\tau^{B \cup C}_N > \max \left(R^B_N,L^C_N \right)$ and notice that it implies that both $R^B_N$ and $L^C_N$ are finite. What we need to show is that for $t > \max \left(R^B_N,L^C_N \right)$ we have $\xi_N (t) \subset \xi^{B \cup C}_N (t)$.

\vspace{0.4 cm}

Let $x \in \xi_N (t)$. There exists $y \in \Z$ such that there is a valid path in the graph $\G$ from $(y,0)$ to $(x,t)$, and we denote it $P_{y \rightarrow x}$. If $y \in B \cup C$ there if of course nothing to prove, so let's assume that $y \not\in B \cup C$. Without loss of generality we assume that $y$ belongs to the right part of the graph $\Z \cap [0,N]$. We consider the right frontier of the left part, denoted $\partial^C_{[-N,\infty[} (t)$, and defined for $t \leq L^C_N$ as follows $$\partial^C_{[-N,\infty[} (t) \text{ } \mydef \text{ } \left\{ (r^C_{[-N,\infty[} (s),s) \in \Z \times \R_+, \text{ } s \in [0,t] \right\},$$ where $r^C_{[-N,\infty[} (t) \text{ } \mydef \text{ } \max \xi^C_{[-N,\infty[}(t)$ for any $t \geq 0$.

\vspace{0.4 cm}

We write simply $\partial^C_{[-N,\infty[}$ for  $\partial^C_{[-N,\infty[} (L^C_N)$. As $L^C_N$ is finite, $\partial^C_{[-N,\infty[}$ goes from $(\max C,0)$ to $(N, L^C)$. Now, as $y > \max C$ and $x \leq N$, $P_{y \rightarrow x}$ need to cross the frontier $\partial^C_{[-N,\infty[}$ at least once, and as in the previous proof, we let $t'$ be the time of the last crossing. Then we can find a valid path $Q_{y' \rightarrow x'}$ for some $y' \in C$, from $(y',0)$, to $(\partial^C_{[-N,\infty[}(t'),t')$, and by concatenation with $P^{t',t}_{y \rightarrow x}$ - the restriction of $P_{y \rightarrow x}$ to the time interval $[t',t]$ - we get a valid path from $(y',0)$ to $(x,t)$, which prove that $x \in \xi^{B \cup C}_N (t)$.

\end{proof}

\begin{lemma} \label{lemma:minmaxfinsemiinf}
For any $t < \tau_N$ we have $$\min \xi_N(t) = \min \xi_{[-N,\infty[} (t),$$ and $$\max \xi_N(t) = \max \xi_{]-\infty, N]} (t).$$
\end{lemma}

\begin{proof}
Let $t < \tau_N$. We prove the lemma only for $\min \xi_N(t)$ as the proof for $\max \xi_N(t)$ is identical. First notice that by monotonicity we have $\min \xi_{[-N,\infty[} (t) \leq \min \xi_N(t)$. Now by definition there exists a path $P_{y \rightarrow x}$ from $(y,0)$ to $(x,t) = (\min \xi_{[-N,\infty[} (t),t)$, for some $y \in \Z \cap [-N,\infty[$. Let $\partial^{\text{left}}_N (t)$ be as in the proof of Lemma \ref{xionlnrn}. If you suppose that $\min \xi_{[-N,\infty[} (t) < \min \xi_N(t)$, then $P_{y \rightarrow x}$ needs to cross  $\partial^{\text{left}}_N (t)$ at least once, so we let $t'$ be the last time of crossing. Then as usual there exists a path $Q_{y' \rightarrow x'}$  for some $y' \in \{-N, \ldots, N \}$, from $(y',0)$, to $(\partial^{\text{left}}_N (t'),t')$, and by concatenation with $P^{t',t}_{y \rightarrow x}$ - the restriction of $P_{y \rightarrow x}$ to the time interval $[t',t]$ - we get a valid path from $(y',0)$ to $(x,t)$, which is a contradiction.
\end{proof}

\vspace{0.4 cm}

\begin{lemma} \label{lemma:A1n}
Let $A \subset \Z$, such that $|A| = n$ for some $n \geq 0$. Then

$$ \P \left(\tau^A < \infty \right) \leq \P \left( \tau^{\{1, \ldots, n\}} < \infty \right).$$
\end{lemma}

\begin{proof}

This lemma is the equivalent for our process to part c of Theorem 1.9 in \cite{ips} concerning the Contact process, and the proof is quite similar so we will merely sketch it here.

\vspace{0.4 cm}

We consider a coupling of $(\xi^{\{1, \ldots n\}}(t))_{t \geq 0}$ with two other processes, denoted $(\xi'(t))_{t \geq 0}$ and $(\eta(t))_{t \geq 0}$. These processes are defined as follows

\begin{itemize}
    \item We have $\xi^{\{1, \ldots n\}}(0) = \{1, \ldots n\}$ by definition, and we let $\xi'(0) = \eta(0) = A$. Active neurons are paired in increasing order in the three processes.
    \item Whenever a leakage occurs in $\xi^{\{1, \ldots n\}}(t)$ at some time $t \geq 0$, the corresponding paired neurons are also affected by leakage  in $\xi'(t)$ and $\eta(t)$.
    \item When a spike occurs in $\xi^{\{1, \ldots n\}}(t)$ at some site $i \in \Z$ for some $t \geq 0$, the spike is propagated in the paired processes in the following sense: the neuron number $i$ becomes quiescent in all three processes ; if neuron $i-1$ (resp. $i+1$) is quiescent in all the processes then it becomes active in all processes, and the newly activated neurons are paired together; if $i-1$ (resp. $i+1$) is already active in $\xi^{\{1, \ldots n\}}(t)$ but is quiescent in $\eta(t)$ then a neuron $i-1$ (resp. $i+1$) becomes active in $\eta(t)$, but not in $\xi'(t)$, and the newly activated neuron then spikes and is affected by leakage effect according to its own independent exponential clocks until further notice ; if $i-1$ (resp. $i+1$) is not already active in $\xi^{\{1, \ldots n\}}(t)$ but is in $\eta(t)$ then neuron $i-1$ (resp. $i+1$) is activated in both $\xi^{\{1, \ldots n\}}(t)$ and $\xi'(t)$, and the newly activated neurons are paired together as well as with the neuron that was already active in $\eta(t)$.
\end{itemize}

With this construction we have that $|\xi^{\{1, \ldots n\}}(t)| = |\xi'(t)| \leq |\eta(t)|$ for any $t \geq 0$, as whenever a neuron is activated in the two first processes, either it is also activated in $\eta(t)$, either the newly activated neuron is paired with an already supernumerary active neuron in $\eta(t)$. Moreover it is clear that $(\eta(t))_{t \geq 0}$ is distributed like $(\xi^A(t))_{t \geq 0}$, and the desired result follows.

\end{proof}

\section{Proof of Theorem \ref{mainth}}

\label{proof}

With this preliminaries completed we can now prove the main result. Note that for reasons related to the way we constructed the proof, even if the theorem is concerned with the time of death of the process $(\xi_N(t))_{t \geq 0}$, the critical value $\gamma'_c$ for which the theorem is stated and here proved is the critical value for the semi-infinite process (see Theorem \ref{theorem:semiinfphasetrans}).

\vspace{0.4 cm}
From the definition of $(\xi_N(t))_{t \geq 0}$ it is clear that, for any $N \in \N$, $\P \left( \tau_N >t \right)$ is continuous and strictly decreasing in $t$ (putting aside the pathological case $\gamma = 0$, in which $\tau_N = \infty$ almost surely), so that we can define $\beta_N$ the unique value in $\R_+$ such that 
$$\P \left( \tau_N > \beta_N \right) = e^{-1}.$$

\vspace{0.4 cm}

We are going to show that

\begin{align}\label{betaNconv}
\frac{\tau_N}{\beta_N} \overset{\mathcal{D}}{\underset{N \rightarrow \infty}{\longrightarrow}} \mathcal{E} (1),
\end{align} from what the result will follow as it will be shown at the end of this proof. The reason for introducing $\beta_N$ is that it will allow us to show that the mean of the exponential random variable is indeed $1$.

\vspace{0.4 cm}

We will prove that the limiting distribution has the memory-less property that characterizes the exponential distribution, that is to say, we will prove that for any $s > 0$ and $t > 0$, we have

\begin{align} \label{result}
\lim_{N \rightarrow \infty} \left| \P\left( \frac{\tau_N}{\beta_N} > s + t \right) - \P \left( \frac{\tau_N}{\beta_N} > s \right)\P\left( \frac{\tau_N}{\beta_N} > t \right) \right| = 0.
\end{align}
  
\vspace{0.4 cm}

So let suppose that $0 < \gamma < \gamma'_c$ and, using the fact that the process is markovian, let us start by observing the following

\begin{align} \label{summationA}
    \nonumber \P \left( \frac{\tau_N}{\beta_N} > s + t \right) &= \sum_{\substack{A \subseteq \{ -N, \ldots N\}\\ A \neq \emptyset}} \P \left( \frac{\tau_N}{\beta_N} > s + t \text{ } \Big| \text{ } \xi_N(\beta_N s) = A \right) \P \Big( \xi_N(\beta_N s) = A \Big)\\
    \nonumber &= \sum_{\substack{A \subseteq \{ -N, \ldots N\}\\ A \neq \emptyset}} \P \left( \frac{\tau_N}{\beta_N} > s + t \text{ } \Big| \text{ } \xi_N(\beta_N s) = A \right) \P \Big( \xi_N(\beta_N s) = A, \text{ }\tau_N > \beta_N s\Big)\\
    \nonumber &= \sum_{\substack{A \subseteq \{ -N, \ldots N\}\\ A \neq \emptyset}} \P \left( \frac{\tau^A_N}{\beta_N} > t\right) \P \Big( \xi_N(\beta_N s) = A, \text{ }\tau_N > \beta_N s\Big)\\
    \nonumber &= \sum_{\substack{A \subseteq \{ -N, \ldots N\}\\ A \neq \emptyset}} \left( \P \left( \frac{\tau^A_N}{\beta_N} > t\right) - \P \left( \frac{\tau_N}{\beta_N} > t\right) \right) \P \Big( \xi_N(\beta_N s) = A, \text{ }\tau_N > \beta_N s\Big)\\
    & \text{ } \text{ } \text{ } +  \P \left( \frac{\tau_N}{\beta_N} > s \right)\P\left( \frac{\tau_N}{\beta_N} > t \right).
\end{align}

\vspace{0.4 cm}
Now, for any $b > 0$, we define the following subset of $\mathcal{P} \left( \Z \right)$

$$ F_b = \left\{A \in \mathcal{P} \left( \Z \right): \frac{|A \cap [-b,0]|}{b+1} > \frac{\rho}{2}, \frac{|A \cap [0,b]|}{b+1} > \frac{\rho}{2} \right\}.$$

Here $\rho$ denotes the density of $(\xi(t))_{t \geq 0}$ as defined earlier. This set is the key point of the proof, the idea behind its definition being that, as the process is spatially ergodic, whenever $b$ will be big enough the measure of the set $F_b$ will be as close to one as needed. 

\vspace{0.4 cm}

Now from (\ref{summationA}), using monotonicity, it follows that

\begin{align} \label{bound}
    \nonumber &\left| \P\left( \frac{\tau_N}{\beta_N} > s + t \right) - \P \left( \frac{\tau_N}{\beta_N} > s \right)\P\left( \frac{\tau_N}{\beta_N} > t \right) \right|\\
    \nonumber &= \sum_{\substack{A \subseteq \{ -N, \ldots N\}\\ A \neq \emptyset}} \left(\P \left( \frac{\tau_N}{\beta_N} > t\right) - \P \left( \frac{\tau^A_N}{\beta_N} > t\right)\right) \P \Big( \xi_N(\beta_N s) = A, \text{ }\tau_N > \beta_N s\Big)\\
    \nonumber &\leq \P \left( \frac{\tau_N}{\beta_N} > s \right)\P\left( \frac{\tau_N}{\beta_N} > t \right) - \sum_{\substack{A \subseteq \{ -N, \ldots N\}\\ A \in F_b}} \P \left( \frac{\tau^A_N}{\beta_N} > t\right) \P \Big( \xi_N(\beta_N s) = A\Big)\\
    \nonumber &\leq \P \left( \frac{\tau_N}{\beta_N} > s \right)\P\left( \frac{\tau_N}{\beta_N} > t \right) - \min_{\substack{A \subseteq \{ -N, \ldots N\}\\ A \in F_b}} \P \left( \frac{\tau^A_N}{\beta_N} > t \right) \P \Big( \xi_N(\beta_N s) \in F_b \Big)\\
    \nonumber &= \P \left( \frac{\tau_N}{\beta_N} > s \right) \Bigg[\P\left( \frac{\tau_N}{\beta_N} > t \right) - \min_{\substack{A \subseteq \{ -N, \ldots N\}\\ A \in F_b}} \P \left( \frac{\tau^A_N}{\beta_N} > t \right)\Bigg]\\
    \nonumber &\text{ } \text{ } \text{ } + \min_{\substack{A \subseteq \{ -N, \ldots N\}\\ A \in F_b}} \P \left( \frac{\tau^A_N}{\beta_N} > t \right) \Bigg[P \left( \frac{\tau_N}{\beta_N} > s \right) -  \P \Big( \xi_N(\beta_N s) \in F_b \Big)\Bigg]\\
    &\leq \Bigg[\P\left( \frac{\tau_N}{\beta_N} > t \right) - \min_{\substack{A \subseteq \{ -N, \ldots N\}\\ A \in F_b}} \P \left( \frac{\tau^A_N}{\beta_N} > t \right)\Bigg] + P \left( \frac{\tau_N}{\beta_N} > s, \text{ } \xi_N(\beta_N s) \not\in F_b \right).
\end{align}

\vspace{0.4 cm}

From now on let fix $\epsilon > 0$. The inequality (\ref{bound}) tells us that in order to prove (\ref{result}) it is sufficient to prove that we can find $b = b_{\epsilon}$ and $N_{\epsilon}$ such that, for all $N \geq N_{\epsilon}$,

\begin{align} \label{outside}
P \Big( \xi_N(\beta_N s) \not\in F_b, \text{ }\tau_N > \beta_N s\Big) < \epsilon,
\end{align} and

\begin{align} \label{inside}
\P\left( \frac{\tau_N}{\beta_N} > t \right) - \min_{\substack{A \subseteq \{ -N, \ldots N\}\\ A \in F_b}} \P \left( \frac{\tau^A_N}{\beta_N} > t \right) < \epsilon.
\end{align}

\vspace{0.4 cm}

We begin with (\ref{inside}). First notice that it is enough to show that there exists $b = b_{\epsilon}$ and $N_{\epsilon}$ such that for all $N \geq N_{\epsilon}$ and all $A \in F_b$ we have

\begin{align*} \P \left( \frac{\tau_N}{\beta_N} > t\right) - \P \left( \frac{\tau^A_N}{\beta_N} > t\right) < \epsilon,
\end{align*} and, using monotonicity again, we have

\begin{align} \label{insidebis}
\P \left( \frac{\tau_N}{\beta_N} > t\right) - \P \left( \frac{\tau^A_N}{\beta_N} > t\right) = \P \left( \frac{\tau_N}{\beta_N} > t, \frac{\tau^A_N}{\beta_N} \leq t\right) \leq \P \Big( \tau_N \neq \tau^A_N\Big),
\end{align} so that it will be sufficient to bound $\P \Big( \tau_N \neq \tau^A_N\Big)$.

\vspace{0.4 cm}

Now for some big enough $n$, as $\tilde{\mu}$ put no mass on $\eta \equiv 0$, we have that

$$ \tilde{\mu}_{[0, +\infty[} \Big( \Big\{B: B \cap [0,n] = \emptyset \Big\} \Big) < \frac{\epsilon}{2}.$$

We take $b_1$ such that $b_1 \cdot \rho / 2 \geq n$ and choose $N_1 > b_1$. Then for any $A \in F_{b_1}$ we have $|A \cap [0,b_1]| \geq b_1 \cdot \rho / 2 \geq n$, so using Lemma \ref{lemma:A1n} we get that for any $A \in F_{b_1}$ and for any $N \geq N_1$

$$ \P \left( \tau^{A\cap[0,b_1]}_{[-N,\infty[} < \infty \right) \leq \P \left( \tau^{\{-N, \ldots ,-N+n\}}_{[-N,\infty[} < \infty \right),$$ and by duality (Theorem \ref{duality}) we have

\begin{align*}
\P \left( \tau^{\{-N, \ldots ,-N+n\}}_{[-N,\infty[} = \infty \right) &= \lim_{t \rightarrow \infty} \P \left( \xi^{\{-N, \ldots ,-N+n\}}_{[-N,\infty[}(t) \neq \emptyset \right)\\
&= \lim_{t \rightarrow \infty} \P \left( C_{[-N,\infty[}(t) \cap \{-N, \ldots ,-N+n\} \neq \emptyset \right)\\
&= \tilde{\mu}_{[0, +\infty[} \left( B: B \cap [0,n] \neq \emptyset \right),
\end{align*} which proves that, for the $b_1$ and $N$ we chose,

\begin{equation} \label{outEleft}
\P \left( \tau^{A\cap[0,b_1]}_{[-N,\infty[} < \infty \right) < \frac{\epsilon}{2}.
\end{equation}

\vspace{0.4 cm}

With the same arguments we also get that

\begin{equation} \label{outEright}
\P \left( \tau^{A \cap [-b_1,0]}_{]-\infty, N]} < \infty \right) < \frac{\epsilon}{2}.
\end{equation}

\vspace{0.4 cm}

This leads us to define the following event

$$ E \text{ } \mydef \text{ } \left\{ \tau^{A\cap[0,b_1]}_{[-N,\infty[} = \infty, \tau^{A\cap[-b_1,0]}_{]-\infty, N]} = \infty \right\}.$$

\vspace{0.4 cm}

We also define the stopping time $$U = \max(L^{A \cap [-b_1,0]}_N, R^{A \cap [0,b_1]}_N),$$ where $L^{A \cap [-b_1,0]}_N$ and $R^{A \cap [0,b_1]}_N$ are as in Lemma \ref{lemma:rl}.

\vspace{0.4 cm}

Then, on $E$ we have $\tau^A_N \geq \tau^{A \cap [-b_1,b_1]}_N \geq U$, and from Lemma \ref{lemma:rl} it follows that, for $t > U$, we have $\xi_N(t) = \xi^{A \cap [-b_1,b_1]}_N (t)$. By monotonicity we have as well for any $t > U$

$$\xi_N(t) = \xi^A_N (t).$$

Therefore, on $E$, we have $\tau_N = \tau_N^A$.

\vspace{0.4 cm}

Now, using this last remark as well as (\ref{outEleft}) and (\ref{outEright}), we get $$\P \Big( \tau_N \neq \tau^A_N \Big) \leq \P \Big( \left\{\tau_N \neq \tau^A_N\right\} \cap E\Big) + \P \left( E^c\right) < \epsilon,$$ which gives a final point to the proof of (\ref{inside}).

\vspace{0.4 cm}

It remains to prove (\ref{outside}). For any choice of $b$, $N$ and $L$ that satisfies the following condition:

\begin{equation}
\label{conditionbnl}
    b < N-L < N,
\end{equation} we have 

\begin{align} \label{outside2}
\nonumber &P \Big( \xi_N(\beta_N s) \not\in F_b, \text{ }\tau_N > \beta_N s\Big)\\
\nonumber &\leq \P \Big( \xi_N(\beta_N s) \not\in F_b, \text{ } \tau_N > \beta_N s, \min \xi_N(\beta_N s) < -N + L, \text{ } \max \xi_N(\beta_N s) > N - L \Big)\\
\nonumber &\text{ } \text{ } \text{ } + \P \Big(\min \xi_N(\beta_N s) \geq -N + L, \text{ } \tau_N > \beta_N s \Big)\\
&\text{ } \text{ } \text{ } + \P \Big(\max \xi_N(\beta_N s) \leq N - L, \text{ } \tau_N > \beta_N s \Big).
\end{align}

\vspace{0.4 cm}

For the first term in the summation in (\ref{outside2}) we have

\begin{align*}
&\P \Big( \xi_N(\beta_N s) \not\in F_b, \text{ } \tau_N > \beta_N s, \text{ } \min \xi_N(\beta_N s) < -N + L, \text{ } \max \xi_N(\beta_N s) > N - L \Big)\\
& \leq \P \Big( \xi_N(\beta_N s) \not\in F_b, \text{ } \tau_N > \beta_N s, \text{ } \min \xi_N(\beta_N s) < -b, \text{ } \max \xi_N(\beta_N s) > b \Big),
\end{align*} so that using Lemma \ref{xionlnrn} this term can be bounded by $\P \Big( \xi(\beta_N s) \not\in F_b\Big)$.

\vspace{0.4 cm}

By the spatial ergodicity of $\mu$ (Theorem \ref{ergodicity}) we have $\mu \left( F^c_b \right) \underset{b \rightarrow \infty}{\longrightarrow} 0$. From this, and using the stochastic monotonicity (Proposition \ref{prop:monotone}), it follows that we can find $b_2$ such that for any $b > b_2$

$$\P \Big( \xi(\beta_N s) \not\in F_{b}\Big) < \frac{\epsilon}{3}.$$

\vspace{0.4 cm}

From Lemma \ref{lemma:minmaxfinsemiinf} we have for any $N$ and $L$

$$\P \Big(\min \xi_N(\beta_N s) \geq -N + L, \text{ } \tau_N > \beta_N s \Big) \leq \P \Big(\min \xi_{[-N,\infty[}(\beta_N s) \geq -N + L \Big).$$

\vspace{0.4 cm}

But by monotone convergence once again we have

\begin{align*}
&\P \Big(\min \xi_{[-N,\infty[}(\beta_N s) \geq -N + L \Big)\\
&\leq \mu_{[-N,+\infty[} \Big( \Big\{ A \subset [-N, \infty[ \text{ } \cap \text{ }\Z : A \cap [-N, -N+L-1] = \emptyset  \Big\}\Big)\\
&= \mu_{[0,+\infty[} \Big( \Big\{ A \subset [0, \infty[ \text{ } \cap \text{ }\Z : A \cap [0, L-1] = \emptyset  \Big\}\Big).
\end{align*}

\vspace{0.4 cm}
And the later will be arbitrarily close to $0$ for arbitrarily big $L$. We therefore fix some big enough $L$ and then fix some $N_2$ such that condition (\ref{conditionbnl}) is satisfied (which in our case means $\max (b_1,b_2) < N_2 - L < L$) and for any  $N \geq N_2$ we get that

$$\P \Big(\min \xi_N(\beta_N s) \geq -N + L, \text{ } \tau_N > \beta_N s \Big) < \frac{\epsilon}{3}.$$

\vspace{0.4 cm}

With the same arguments and by symmetry the last term in (\ref{outside2}) is also bounded by $\epsilon/3$ for $N \geq N_2$ and for the same choice of $L$. To finish of course we take $b_{\epsilon} = \max (b_1,b_2)$ and $N_\epsilon = \max(N_1,N_2)$ and both (\ref{inside}) and (\ref{outside}) are satisfied.

\vspace{0.4 cm}

Notice that formally, from (\ref{result}) alone, there is actually three possibilities for the survival function of the limiting distribution: it could be identically equal to $1$, identically equal to $0$ or it could be $t \mapsto e^{- t}$. Nonetheless the definition of $\beta_N$ rules out the two first cases, which correspond to a degenerate random variable, so that we indeed have

$$ \lim_{N \rightarrow \infty} \P\left( \frac{\tau_N}{\beta_N} > t \right) = e^{-t}.$$

\vspace{0.4 cm}

It only remains to show that $\beta_N$ can be replaced by $\E \left( \tau_N \right)$. We know from (\ref{summationA}) and by monotonicity that for any $s, t \geq 0$ we have $$ \P\left( \frac{\tau_N}{\beta_N} > s + t \right) \leq \P \left( \frac{\tau_N}{\beta_N} > s \right)\P\left( \frac{\tau_N}{\beta_N} > t \right),$$ so that it follows from the definition of $\beta_N$ that for any integer $n$ we have

$$ \P\left( \frac{\tau_N}{\beta_N} > n \right) \leq e^{-n}.$$

\vspace{0.4 cm}

In general for any $t \geq 0$ we therefore have

$$ \P\left( \frac{\tau_N}{\beta_N} > t \right) \leq e^{- \lfloor t \rfloor }.$$

\vspace{0.4 cm}

Moreover we have

$$ \frac{\E \left( \tau_N \right)}{\beta_N} = \int_0^{\infty} \P \left(\frac{\tau_N}{\beta_N} > t \right) dt,$$ so that finally - by Dominated Convergence Theorem - we get

$$ \lim_{N \rightarrow \infty} \frac{\E \left( \tau_N \right)}{\beta_N} = \int_0^{\infty} \lim_{N \rightarrow \infty} \P \left(\frac{\tau_N}{\beta_N} > t \right) dt = \int_0^{\infty} e^{-t} dt = 1.$$

\section{Annex}

In this annex we prove that the set of continuous and increasing function on $\{0,1\}^\Z$ is distribution determining, as it was claimed in Section \ref{sectionasinfproc}.

\begin{proposition}
Let $\mu_1$ and $\mu_2$ be two probability measures on $(\{0,1\}^\Z, \mathcal{B}(\{0,1\}^\Z))$, where $\mathcal{B}(\{0,1\}^\Z)$ denotes the borelian $\sigma$-algebra on $\{0,1\}^\Z$. Suppose that for any continuous and increasing function $f: \{0,1\}^\Z \mapsto \R$ we have the following equality

$$\int f d \mu_1 = \int f d \mu_2,$$

then $\mu_1 = \mu_2$.
\end{proposition}

\begin{proof}
Fix an arbitrary $B \in \mathcal{B}(\{0,1\}^\Z)$ different from the empty set. For any $n \geq 0$ define the function $f_n$ on $\{0,1\}^\Z$ by

$$ f_n(x) = \mathbbm{1}\{x_{[-n,n]} \geq b_{[-n,n]} \text{ for some } b \in B\},$$

where for any $x \in \{0,1\}^\Z$, $x_{[-n,n]}$ is defined as follows

\[x_{[-n,n]}(i)= 
\begin{cases}
    x(i)& \text{if } i \in [-n,n],\\
    0 & \text{otherwise}.
\end{cases}
\]

\vspace{0.4 cm}

The function $f_n$ can be written as follows

\begin{align*}
    f_n(x) \text{ } &= \text{ } g_n(x) + h_n(x)\\
    &\mydef \text{ }\mathbbm{1}\{x_{[-n,n]} = b_{[-n,n]} \text{ for some } b \in B\} + \mathbbm{1}\{x_{[-n,n]} > b_{[-n,n]} \text{ for some } b \in B\}.\\
\end{align*}

\vspace{0.4 cm}

Notice that $f_n$ is a continuous and increasing function, so that by hypothesis we have $$ \int f_n d \mu_1 = \int f_n d \mu_2,$$ which can be written $$ \int g_n d \mu_1 + \int h_n d \mu_1 = \int g_n d \mu_2 + \int h_n d \mu_2.$$

\vspace{0.4 cm}

Now $h_n$ is a continuous and increasing function as well, so that the two integrals involving $h_n$ simply cancel out in the previous equation, and we are left with

$$ \int g_n d \mu_1 = \int g_n d \mu_2.$$

\vspace{0.4 cm}

Finally we have $\lim_{n \rightarrow \infty} g_n (x) =  \mathbbm{1}\{x \in B\}$, so that by dominated convergence theorem $$ \mu_1(B) = \lim_{n \rightarrow \infty} \int g_n d \mu_1 = \lim_{n \rightarrow \infty} \int g_n d \mu_2 = \mu_2(B),$$ which ends the proof.

\end{proof}

\newpage

\section*{Acknowledgements}
This work is part of my PhD thesis. I thank my PhD adviser Antonio Galves for introducing me to the subject of metastability and for fruitful discussions. Many thanks also to Christophe Pouzat who introduced me to the field of neuromathematics. This article was produced as part of the activities of FAPESP  Research, Innovation and Dissemination Center for Neuromathematics (grant number 2013/07699-0 , S.Paulo Research Foundation), and the author was supported by a FAPESP scholarship (grant number 2017/02035-7).

\end{document}